\documentclass{amsart}
\usepackage{amsmath, amssymb, mathrsfs, verbatim, multirow}
\usepackage[all]{xy}

\newtheorem{Teo}{Theorem}[section]

\theoremstyle{definition}
\newtheorem{Def}[Teo]{Definition}

\newtheorem{Obs}[Teo]{Remark}

\newcommand{\R}{\mathbb{R}}

\newcommand{\N}{\mathbb{N}}

\newcommand{\lra}{\longrightarrow}

\newcommand{\VR}{\mathcal{O}}
\newcommand{\PI}{\mathfrak{p}}
\newcommand{\MI}{\mathfrak{m}}

\newcommand{\hei}{\mbox{\rm ht}}
\newcommand{\QF}{\mbox{\rm Quot}}

\newcommand{\SU}{\mbox{\rm supp}}

\newcommand{\ann}{\mbox{\rm ann}}

\newcommand{\rk}{\mbox{\rm rk}}
\newcommand{\Sp}{\mbox{\rm Spec}}

\begin{document}
\title{On the local uniformization problem}
\author{Josnei Novacoski}
\address{Josnei Novacoski\newline\indent CAPES Foundation  \newline \indent Ministry of Education of Brazil \newline \indent Bras\'ilia/DF 70040-020 \newline \indent Brazil}
\email{jan328@mail.usask.ca}

\author{Mark Spivakovsky}
\address{Mark Spivakovsky\newline\indent
CNRS and Institut de Math\'ematiques de Toulouse\newline\indent
Universit\'e Paul Sabatier\newline\indent
118, route de Narbonne\newline\indent
31062 Toulouse cedex 9\newline\indent
France}
\email{mark.spivakovsky@math.univ-toulouse.fr}

\thanks{During the realization of this project the first author was supported by a research grant from the program ``Ci\^encia sem Fronteiras" from the Brazilian government.}

\begin{abstract}
In this paper we give a short introduction to the local uniformization problem. This follows a similar line as the one presented by the second author in his talk at ALANT 3. We also discuss our paper on the reduction of local uniformization to the rank one case. In that paper, we prove that in order to obtain local uniformization for valuations centered at objects of a subcategory of the category of noetherian integral domains, it is enough to prove it for rank one valuations centered at objects of the same category. We also announce an extension of this work which was partially developed during ALANT 3. This extension says that the reduction mentioned above also works for noetherian rings with zero divisors (including the case of non-reduced rings).
\end{abstract}

\maketitle

\section{Introduction}
Resolution of singularities for an algebraic variety is an important branch of algebraic geometry. Roughly speaking, for an algebraic variety $V$ a resolution of singularities is an algebraic variety $V'$ with no singularities, birationally equivalent to $V$ (a precise definition is provided in Section \ref{resolunifo}). Local uniformization is the local version of resolution of singularities. Namely, given a valuation $\nu$ centered at a point $\PI\in V$ we want to find $V'$ birationally equivalent to $V$ such that the center $\PI'$ of $V'$ is non-singular (see Section \ref{resolunifo}).

Local uniformization was introduced by Zariski in order to prove resolution of singularities. His approach consists of two steps: proving local uniformization for every valuation and use these local solutions to obtain a resolution of all singularities. In this second step, the quasi-compactness of the Zariski topology on the space of valuations plays an important role. This is because when we obtain local uniformization for a given valuation, every valuation in an open neighbourhood of it is resolved. This can be seen as follows. For an irreducible algebraic variety $V$, the set $U(V)$ of valuations of $K(V)$ having a center in $V$ is an open set in the Zariski topology. Moreover, the map $g_V:U(V)\lra V$, taking a valuation to its center is continuous by definition of Zariski topology on the Riemann-Zariski space of valuations. If the center of a given valuation $\nu$ in $V$ is non-singular, then the pre-image of the set of non-singular points of $V$ (which is an open set) is an open neighbourhood $U_\nu$ of $\nu$ such that every valuation in $U_\nu$ has a non-singular center in $V$. Hence, using the quasi-compactness of the Zariski space of valuations, it is enough to ``glue" only finitely many solutions.

Zariski succeeded in 1940 (see \cite{Zar}) in proving local uniformization for valuations having a center on any algebraic variety over a field of characteristic zero. He used this to prove resolution of singularities for algebraic varieties of dimension smaller or equal to three over a field of characteristic zero. Abhyankar proved in 1956 (see \cite{Ab_4}), using Zariski's approach, resolution of singularities for algebraic varieties of dimension 3 over a perfect field of characteristic $p\geq 6$. In the last several years, Cossart and Piltant concluded the proof (see \cite{Cos1}, \cite{Cos2}, \cite{Cos3}) of resolution of singularities for any
quasi-excellent noetherian scheme of dimension 3 (they also used Zariski's approach).

Resolution of singularities is known for algebraic varieties over a field of characteristic zero. The first full proof of it was given in \cite{Hir_1} by Hironaka (for this work he received a Fields Medal in 1970). However, both resolution of singularities and local uniformization are open problems for algebraic varieties of dimension greater than 3 and positive characteristic.

In the cases where local uniformization is known the proof usually proceeds by induction on the rank of the valuation. In \cite{Nov} we proved that this process does not depend on the proof for the rank one case. Namely, we proved that if every rank one valuation centered in some ring of a given category admits local uniformization, then every valuation centered at objects of this category admits local uniformization (see Section \ref{Redlunif} for a more precise statement). In \cite{Nov} we only dealt with valuations centered at integral domains. However, in our recently submitted paper \cite{Nov1}, we extend this result to valuations centered at any type of noetherian local rings (for instance, such rings might have nilpotent elements).

In Section \ref{resolunifo}, we present some basic definitions from algebraic geometry. We give the background necessary in order to define resolution of singularities and local uniformization. For simplicity of exposition, in this introduction we chose to restrict ourselves to the special case of algebraic varieties over an algebraically closed field. A reader familiar with scheme theory will recognize the fact that all the concepts and definitions, as well as our results on reduction to the rank one case, generalize naturally to a much more general setting of arbitrary noetherian schemes.

In Section \ref{Redlunif} we describe our results on the reduction of local uniformization to the rank one case.

\par\medskip
\textbf{Acknowledgements.} We thank the anonymous referee for a careful reading of the paper and for providing many useful suggestions and corrections.

\section[Singularities and uniformization]{Resolution of singularities and local uniformization}\label{resolunifo}

Take an algebraically closed field $K$ and a prime ideal $I=(f_1,\ldots,f_r)$ of $K[\underline X]:=K[X_1,\ldots,X_n]$. The
\textbf{affine algebraic variety $V(I)$} is defined as the set of zeros of $I$ in $K^n$:
\[
V(I):=\{p\in K^n\mid f_i(p)=0\mbox{ for every }i, 1\leq i\leq r\}.
\]
We define the \textbf{coordinate ring of $V:=V(I)$} as the ring $K[V]:=K[X_1,\ldots,X_n]/I$. The \textbf{function field $K(V)$ of $V$} is defined as the quotient field of $K[V]$. Set $x_i=X_i+I$ for each $i=1,\ldots, n$. Then $K[V]=K[\underline x]:=K[x_1,\ldots,x_n]$ and
$K(V)=K(\underline x):=K(x_1,\ldots,x_n)$.

We define a topology $\mathcal Z(V)$ on $V$ by setting as open sets the sets of the form $V\setminus V(\mathcal I)$ where $\mathcal I$ runs over all the ideals of $K[\underline X]$. A function $f:V\lra K$ is said to be \textbf{regular} if there exist $U_f\in \mathcal
Z(V)$ and $g,h\in K[V]$ such that $h(p)\neq 0$ and $f(p)=g(p)/h(p)$ for every $p\in U_f$. Then every element in $K(V)$ can be seen as the equivalence class of a pair $(U_f,f)$ under the equivalence given by $(U_f,f)\sim (U_{f'},f')$ if and only if $f$ coincides with
$f'$ in $U_f\cap U_{f'}$. A continuous function $\Phi:V\lra V'$ between the varieties $V$ and $V'$ (over $K$) is a \textbf{morphism} if for every regular function $f: V'\lra K$ the function $f\circ \Phi: \Phi^{-1}(V')\lra K$ is also regular. Let $\mathcal M(V,V')$ denote the set of morphisms from $V$ to $V'$.

A \textbf{rational map} is an equivalence class of pairs of the form
$$
(U,\Phi)\in \mathcal Z(V)\times \mathcal M(V,V')
$$
under the equivalence $(U,\Phi)\sim (U',\Phi')$ if and only if $\Phi$ and $\Phi'$ coincide in $U\cap U'$. Hence, every rational map
$V\lra V'$ induces a $K$-homomorphism $K(V')\lra K(V)$ of fields. A rational map is said to be \textbf{birational} if the induced map on the function fields is an isomorphism.

For a point $p\in V$ we consider the ideal $\PI=\{f(\underline x)\in K[V]\mid f(p)\neq 0\}$ of $K[V]$. This is a prime ideal of $K(V)$. Define the local ring of $V$ at $p$ as
\[
\VR_{V,p}:=\{\phi\in K(V)\mid \phi=f(x)/g(x)\mbox{ such that }g(p)\neq 0\}=K[V]_\PI.
\]
The point $p$ is said to be \textbf{regular} (or \textbf{non-singular}) if $\VR_{V,p}$ is regular, i.e., if the only maximal ideal of $\VR_{V,p}$ is generated by $\dim(\VR_{V,p})$-many elements.

An algebraic variety $V$ over $K$ is said to admit \textbf{resolution of singularities} if there exists a proper birational morphism
$\pi:V'\lra V$ such that every point of $V'$ is regular. Proper means that for every valuation $\nu$ of $K(V)=K(V')$ if $\VR_{V,p}$ is dominated by $\VR_\nu$, then there exists a unique $p'\in\pi^{-1}(p)$ such that $\VR_{V',p'}$ is dominated by $\VR_\nu$. If $\nu$ is a valuation of $K(V)$, then the pair $(V,\nu)$ is said to admit \textbf{local uniformization} if there exists a proper birational morphism $\pi:V'\lra V$ from a variety $V'$ to $V$ such that the unique point $p'\in\pi^{-1}(p)$ for which
$\VR_{V',p'}\subseteq\VR_\nu$ is regular.

In the modern language, an algebraic variety over $K$ is replaced by an ``integral separated scheme of finite type over $K$". A scheme is the analogue in algebraic geometry of a	manifold. A manifold is a topological space which is the union of open sets which are homeomorphic to $\R^n$. Analogously, a scheme of finite type over $K$ with function field $F$ is the union of objects which are isomorphic (in the category of schemes) to spaces (called affine schemes) of the form $\Sp(A)$ where
\[
A=K[x_1,\ldots,x_n]\subseteq F\mbox{ such that }F=K(x_1,\ldots,x_n).
\]
We can think of $A$ as the coordinate ring of an algebraic variety $V$. Since $\Sp(A)$ consists of all the prime ideals of $A$, we are extending the definition of point to include also irreducible subvarieties of $V$. Since local uniformization is a local problem, it is enough to consider only affine varieties. Also, it is natural to consider valuations whose center at such a variety may be the generic point of an irreducible subvariety (rather than only a ``classical'' point).

A natural way of resolving singularities is by blowing ups. Take $V=\Sp(A)$ and consider a valuation $\nu$ of $F:=\QF(A)$ having a center $\PI'$ on $V$ (i.e., $A\subseteq \VR_\nu$ and $\PI':=A\cap \MI_\nu$). We can describe a blowing up of $V$ at $\PI'$ along $\nu$ in the following way. Let $R$ be a noetherian local domain (e.g., $R=A_{\PI'}$) and a valuation $\nu$ of $\QF(R)$ centered at $R$
(i.e., such that $R$ is dominated by $\VR_\nu$). A \textbf{local blowing up of $R$ with respect to $\nu$} is an inclusion map $R\lra
R^{(1)}$ where
\[
R^{(1)}:=R[a_1,\cdots, a_l]_{\MI_\nu\cap R[a_1,\cdots, a_l]}\mbox{ for some }a_1,\ldots,a_l\in \VR_\nu.
\]
Let $X$ and $X'$ be two schemes and consider a proper birational map $\Phi:X'\lra X$. Let $\nu$ be a valuation of $K(X)=K(X')$ having a center $x$ in $X$. Let $x'\in X'$ be the center of $\nu$ in $X'$ (which exists by the valuative criterion for properness). Let $V=\Sp(A)\subseteq X$ and $V'=\Sp(A')\subseteq X'$ be open subsets such that $x\in V$ and $x'\in V'$. Then the induced map
$\VR_{X,x}\lra\VR_{X',x'}$ is a local blowing up. Hence, we modify the definition of local uniformization in this more modern language to the following: for a noetherian local domain $R$ and a valuation $\nu$ of $F=\QF(R)$ centered at $R$, we say that the pair $(R,\nu)$ admits local uniformization if there exists a local blowing up $R\lra R^{(1)}$ of $R$ with respect to $\nu$ such that $R^{(1)}$ is regular.

\section[Reduction]{Reduction of local uniformization to the rank one case}\label{Redlunif}

An important step in most of the proofs of local uniformization is to reduce local uniformization to the case of rank one valuations, i.e., valuations whose value group can be embedded in $\mathbb R$ as an ordered group. This step appears, for instance, in Zariski's proof for the characteristic zero case and in Cossart and Piltant's proof for positive characteristic and dimension at most $3$. Let us formulate this reduction more precisely.

Let $\mathcal N$ denote the category of all noetherian local domains. Let $\mathcal N'$ be a subcategory of $\mathcal N$. We will say that $\mathcal N'$ admits \textbf{reduction to rank one} if the following holds: if every rank one valuation centered at any object of
$\mathcal N'$ admits local uniformization, then all the valuations centered at objects of $\mathcal N'$ admit local uniformization. In \cite{Nov}, we prove the following:
\begin{Teo}\label{mainsthmspinov1}
Let $\mathcal{N'}$ be a subcategory of $\mathcal{N}$ which is closed under taking homomorphic images and localizing at a prime ideal any finitely generated birational extension. Then $\mathcal N'$ admits reduction to rank one.
\end{Teo}

We sketch the proof of Theorem \ref{mainsthmspinov1} here for the convenience of the reader. For a complete proof of it, see
\cite{Nov}. We use induction on the rank of the valuations. Take $n\in\N$. If $n=1$, then by assumption all the valuations of rank $n$ centered at members of $\mathcal{N'}$ admit local uniformization. If $n>1$, then we assume that all the valuations of rank smaller than $n$ centered at members of $\mathcal{N'}$ admit local uniformization. Then we prove that also valuations of rank $n$ centered at members of $\mathcal{N'}$ admit local uniformization.

Fix a valuation of rank $n$ centered at some member $R$ of $\mathcal{N'}$. Decompose $\nu$ as $\nu=\nu_1\circ\nu_2$ such that $\nu_1$ and $\nu_2$ have rank smaller than $n$. We write $\PI$ for the center of $\nu_1$ in $R$. Then $\nu_1$ is a valuation centered at $R_\PI$ and $\nu_2$ is a valuation of $R_\PI/\PI R_\PI$ centered at $R/\PI$. By the induction hypothesis we can assume that $\nu_1$ and $\nu_2$ admit local uniformization. This implies that there exist local blowing ups along $\nu_1$ and $\nu_2$ such that, after these local blowing ups, the rings $R_\PI$ and $R/\PI$ become regular. The first main step is to prove that these local blowing ups can be lifted to local blowing ups along $\nu$. For more precise statements, as well as their proofs, see Corollaries 2.17 and 2.20 of \cite{Nov}.

By the discussion in the previous paragraph we can assume (replacing $R$ by some blowing up $R^{(1)}$ of $R$) that $R_\PI$ and $R/\PI$ are regular. Let $(y_1,\ldots,y_r)\subseteq\PI$ be a regular system of parameters for $\mathfrak{p}R_\PI$ and $x_1,\ldots, x_t$ a set of elements of $R\setminus\mathfrak{p}$, whose images modulo $\mathfrak{p}$ form a regular system of parameters of $\mathfrak{m}/\mathfrak{p}$. If $y_1,\ldots,y_r$ generate $\PI$, then $R$ is regular. Indeed, since $y_1,\ldots,y_r,x_1,\ldots,x_t$ generate $\MI$ we have $r+t\geq \dim (R)$. Also, since $r= \dim \left(R_\mathfrak{p}\right)=\hei
\left(\PI\right)$ and $t=\dim \left(R/\PI\right)=\hei\left(\MI/\PI\right)$ we have
\[
\dim (R)=\hei\left(\MI\right)\geq \hei\left(\PI\right)+\hei\left(\MI/\PI\right)=r+t\geq \dim (R).
\]
Therefore, $r+t=\dim (R)$ and $y_1,\ldots,y_r,x_1,\ldots,x_t$ is a minimal set of generators of $\MI$, hence $\left(R,\MI\right)$ is regular.

If $y_1,\ldots,y_r$ do not generate $\PI$, take $y_{r+1},\ldots,y_{r+s}\in \PI$ such that $y_1,\ldots,y_r$, $y_{r+1},\ldots,y_{r+s}$ generate $\PI$. The goal is to eliminate these elements in order to achieve the situation in the previous paragraph. We will sketch the idea used in the proof by eliminating $y_{r+1}$. Since the residues of $y_1,\ldots, y_r$ modulo $\left(\PI R_\PI\right)^2$ form an $R_\PI/\PI R_\PI$-basis of $\PI R_\PI/\left(\PI R_\PI\right)^2$ we can find an equation
\[
a_1y_{r+1}+b_{11}y_1+\ldots +b_{r1}y_r-h_1=0
\]
where $a_1\in R\setminus \PI$ and $h_1\in  \left(y_1,\ldots ,y_r\right)^2$. Then we blow up $R$ with respect to $\nu$ along the ideal $\left(a_1,y_1,\ldots,y_r\right)$ obtaining a new local domain $R^{(1)}$. In $R^{(1)}$ the new regular system of parameters
$\left(y^{(1)}_1,\ldots, y^{(1)}_r\right)$ of $R^{(1)}_{\PI^{(1)}}$ is obtained by
$$
y_1=a_1y_1^{\left(1\right)}, y_2=a_1y_2^{\left(1\right)},\ldots, y_r=a_1y_r^{\left(1\right)}.
$$
We can rewrite the previous relation as
\begin{equation}\label{eqthatshimpdomain}
a_1\left(y_{r+1}+b_{11}y_1^{\left(1\right)}+\ldots+b_{r1}y_r^{\left(1\right)}-h_1'\right)=0
\end{equation}
where $h_1=a_1h_1'$ with $h_1'\in \left(y_1^{(1)},\ldots,y_r^{(1)}\right)^2$. In particular,
\begin{equation}\label{eqobtainedafter}
y_{r+1}+b_{11}y_1^{\left(1\right)}+\ldots+b_{r1}y_r^{\left(1\right)}-h_1'=0,
\end{equation}
which implies that $y_{r+1}\in \left(y^{(1)}_1,\ldots,y^{(1)}_r\right)$.

\begin{Obs}\label{rmkondiffic1}
In equation (\ref{eqthatshimpdomain}), we used the fact that $R$ is a domain to obtain (since $a_1\neq 0$) equation (\ref{eqobtainedafter}). This is the first difficulty to overcome when ``adapting" our proof to rings which are not necessarily domains. 
\end{Obs}
In \cite{Nov}, we also consider stronger forms of local uniformization. Two of those are what we called \textbf{weak embedded local uniformization} and \textbf{embedded local uniformization}. Weak embedded local uniformization asks whether for every given finite subset $\mathcal F$ of $\VR_\nu$ we can find a local blowing up as in $R\lra R^{(1)}$ such that $R^{(1)}$ is regular and a regular system of parameters $u=(u_1,\ldots,u_d)$ of $R^{(1)}$ such that all elements of $\mathcal F$ are monomials in $u$. We order the elements of the set $\mathcal F$ above by their values, i.e., $\mathcal F=\{f_1,\ldots,f_q\}$ such that $\nu(f_1)\leq\ldots\leq\nu(f_q)$. We can ask whether we can find a regular local domain $R^{(1)}$ with regular system of parameters $u$ as before such that the elements $f_i$ are monomials in $u$ and moreover, $f_1\mid_{R^{(1)}}\ldots\mid_{R^{(1)}} f_q$. This version is called embedded local uniformization. Recently, it was noted to us by Schoutens (in an e-mail by Kuhlmann), that these two concepts are equivalent.

We present a sketch of the proof of this equivalence. Embedded local uniformization implies, from its definition, weak embedded local uniformization. For the converse, assume that weak embedded local uniformization is true for the pair $(R,\nu)$. Take a subset $\mathcal F=\{f_1,\ldots, f_q\}$ of $\VR_\nu$. Define
\[
\mathcal F'=\mathcal F\cup \{f_i-f_j\mid 1\leq j< i\leq q\}.
\]
By assumption, there exists a local blowing up $R\lra R^{(1)}$ such that $R^{(1)}$ is regular with a regular system of parameters $u=(u_1,\ldots, u_d)$ such that every element in $\mathcal F'$ is a monomial in $u$. For fixed $i$ and $j$, $1\leq j<i\leq q$, we have that $f_i=au^{\alpha}$, $f_j=bu^{\beta}$ and $f_i-f_j=cu^\gamma$ for some units $a,b,c\in R^{(1)}$ and $\alpha,\beta,\gamma\in (\N\cup\{0\})^d$. This means that
\[
cu^\gamma=f_i-f_j=au^\alpha-bu^\beta=\left(au^{\alpha-\sigma}-bu^{\beta-\sigma}\right)u^\sigma,
\]
where $\sigma=(\min\{\alpha_1,\beta_1\},\ldots,\min\{\alpha_d,\beta_d\})$. This is only possible if $\sigma=\gamma$ and $au^{\alpha-\sigma}-bu^{\beta-\sigma}$ is a unit in $R^{(1)}$, which implies that $\beta=\sigma$. Therefore, $f_j$ divides $f_i$ in $R^{(1)}$, which is what we wanted to prove.

During ALANT 3 we have made important developments in our now submitted paper \cite{Nov1}. The aim of this paper is to generalize Theorem \ref{mainsthmspinov1} to the case where objects of $\mathcal N$ are not necessarily integral domains. For that case we have to adapt our definitions.

Take a noetherian local ring $R$ ($R$ may have zero divisors and even nilpotent elements) and an abelian group $\Gamma$. Take $\infty$ to be an element not in $\Gamma$ and set $\Gamma_\infty$ to be $\Gamma\cup\{\infty\}$ with extensions of addition and order as usual.

\begin{Def}
A \textbf{valuation on $R$} is a map $\nu:R\lra \Gamma_\infty$ such that the following holds:
\begin{description}
\item[(V1)] $\nu(ab)=\nu(a)+\nu(b)$ for every $a,b\in R$,
\item[(V2)] $\nu(a+b)\geq \min\{\nu(a),\nu(b)\}$ for every $a,b\in R$,
\item[(V3)] $\nu(1)=0$ and $\nu(0)=\infty$,
\item[(V4)] $\SU(\nu):=\{a\in R\mid \nu(a)=\infty\}$ is a minimal prime ideal of $R$.
\end{description}
\end{Def}

We observe that if $S$ is a multiplicative set of $R$, contained in $R\setminus \SU(\nu)$, then $\nu$ extends to a valuation on $R_S$
(which we denote again by $\nu$) by setting $\nu(a/b)=\nu(a)-\nu(b)$. A valuation $\nu$ on $R$ is said to have a center if $\nu(a)\geq 0$ for every $a\in R$. In this case the center is defined as $\mathfrak c_{R}(\nu):=\{a\in R\mid\nu(a)>0\}$. Moreover, if $R$ is a local ring with unique maximal ideal $\MI$, then a valuation $\nu$ on $R$ is said to be centered at $R$ if $\nu(a)\geq 0$ for every $a\in R$ and $\nu(a)>0$ for every $a\in \MI$. We observe that if $\nu$ is a valuation having a center at $R$, then the extension of $\nu$ to
$R_{\mathfrak c_R(\nu)}$ is centered at $R_{\mathfrak c_R(\nu)}$.

Take an element $b\in R\setminus \SU(\nu)$. If we consider the natural map $\Phi:R\lra R_b$ given by $\Phi(a)=a/1$, then
\[
J(b):=\ker \Phi=\bigcup_{i=1}^\infty \ann_R(b^i).
\]
Hence, we have a natural embedding $R/J(b)\subseteq R_{b}$. Take $a_1,\ldots,a_r\in R$ such that $\nu(a_i)\geq \nu(b)$ for each $i$,
$1\leq i\leq r$. Consider the subring $R':=R/J(b)[a_1/b,\ldots,a_r/b]$ of $R_b$. Then $\nu$ has a center $\mathfrak c_{R'}(\nu)$ on
$R'$. Consider the ring $R^{(1)}:=R'_{\mathfrak c_{R'}(\nu')}$ and the extension $\nu^{(1)}$ of $\nu$ to $R^{(1)}$.

\begin{Def}
A local blowing up of $R$ along $\nu$ is an inclusion map $R\lra R^{(1)}$ of the form described above.
\end{Def}
Let $I=\sqrt{(0)}$ ($I$ is the \textbf{nilradical of $R$}).
\begin{Obs} If $b\in I$ then $R_b$ is the zero ring (that is, the one-element ring in which 0=1). If $b\in\SU(\nu)\setminus I$ then the ring $R_b$ and the homomorphism $\Phi:R\rightarrow R_b$ are well defined but there does not exist a localization $R^{(1)}$ of $R_b$ and a valuation $\nu^{(1)}$ centered in $R^{(1)}$ whose restriction to $R$ is $\nu$. This is why we limit ourselves to the case
$b\notin\SU(\nu)$.
\end{Obs}
\begin{Obs} The ring $R_{\SU(\nu)}$ has only one associated prime ideal. If $R_{\SU(\nu)}$ contains non-zero nilpotent elements, then so does every local blowing up $R^{(1)}$ of $R$. Therefore in this case there is no hope of making $R^{(1)}$ regular; the best we can ask is for $\left(R^{(1)}\right)_{\text{red}}$ to be regular. Furthermore, it is both natural and possible to look  or some form of ``constant'' or ``uniform'' behaviour of the nilradical of $R^{(1)}$ along $\Sp \left(R^{(1)}\right)_{\text{red}}$. The intuitive idea of uniform behaviour of a module along a scheme in algebraic geometry is often embodied in the concept of flatness. Therefore, in order to define local uniformization for rings with nilpotents it is natural to ask that the nilradical of $R^{(1)}$ become a flat $\left(R^{(1)}\right)_{\text{red}}$-module. In fact, one can do slightly better and require not only the nilradical itself but all of the successive quotients of its powers to be flat. A finitely generated module over a local ring is flat if and only if it is free. These considerations motivate the following definition.
\end{Obs}
\begin{Def} Assume that $R_{\rm red}$ is regular. We say that $\Sp(R)$ is \textbf{normally flat along $\Sp(R_{\rm red})$} if
$I^n/I^{n+1}$ is an $R_{\rm red}$-free module for every $n\in\N$.
\end{Def}
Since $I^n=(0)$ for every $n>N$ for some $N\in\N$, this condition is equivalent to the freeness of the finitely many modules
$I/I^2,\ldots, I^N/I^{N+1}=I^N$.

A pair $(R,\nu)$ as above is said to admit \textbf{local uniformization} if there exists a local blowing up $R\lra R^{(1)}$ along $\nu$ such that $R^{(1)}_{\rm red}$ is regular and $\Sp(R^{(1)})$ is normally flat along $\Sp\left(R^{(1)}\right)_{\rm red}$. 

Consider the category $\mathcal M$ of all noetherian rings (not necessarily integral domains). Again, we will say that $\mathcal M'$ admits \textbf{reduction to rank one} if the following holds: if every rank one valuation centered at any member of $\mathcal M'$ admits local uniformization, then all the valuations centered at members of $\mathcal M'$ admit local uniformization. The following is the main theorem of \cite{Nov1}.

\begin{Teo}\label{thmredcasenilp}
Let $\mathcal{M'}$ be a subcategory of $\mathcal{M}$ which is closed under taking homomorphic images and localizing at a prime any finitely generated birational extension. Then $\mathcal M'$ admits reduction to rank one.
\end{Teo}

The proof of Theorem \ref{thmredcasenilp} consists of three main steps. The first step is to prove that for every local ring $R$ and every valuation $\nu$ centered on $R$, there exists a local blowing up $R\lra R^{(1)}$ such that $R^{(1)}$ has only one associated prime ideal. This first step eliminates the difficulty described in Remark \ref{rmkondiffic1}.

Then we consider a decomposition $\nu=\nu_1\circ\nu_2$ of $\nu$ such that $\rk(\nu_1)<\rk(\nu)$ and $\rk(\nu_2)<\rk(\nu)$. Using induction, we can assume that both $\nu_1$ and $\nu_2$ admit local uniformization. The second main step consists in using this to prove that there exists a local blowing up $R^{(1)}\lra R^{(2)}$ such that $\left(R^{(2)}\right)_{\rm red}$ is regular. In this step we use a similar procedure as the one used to prove Theorem \ref{mainsthmspinov1}.

The third and final step is to prove that there exists a further local blowing up $R^{(2)}\lra R^{(3)}$ such that $\left(R^{(3)}\right)_{\rm red}$ is regular and $I_{(3)}^n/I_{(3)}^{n+1}$ is an $\left(R^{(3)}\right)_{\rm red}$-free module for every $n\in\N$ (here $I_{(3)}$ denotes the nilradical of $R^{(3)}$). The method used here is again similar to the one used in the proof of Theorem \ref{mainsthmspinov1}. The main point is that if a module $I_{(2)}^n/I_{(2)}^{n+1}$ is not free, then we can obtain some equation similar to equation (\ref{eqthatshimpdomain}). This equation gives us the ideal along which we have to blow up.

In the second and third steps described above we have to be careful when choosing the ideals because we do not want to destroy the good properties already achieved. This is the main difficulty in proving Theorem \ref{thmredcasenilp}. This is why, for instance, we cannot reduce (part of) its proof to Theorem \ref{mainsthmspinov1}. For a detailed proof of Theorem \ref{thmredcasenilp}, see
\cite{Nov1}.

\end{document}